\documentclass[10pt,a4paper]{article}
\usepackage[T1]{fontenc}
\usepackage[utf8]{inputenc}
\usepackage{mathptmx}
\usepackage[margin=25mm]{geometry}
\usepackage{amsmath,amssymb}
\usepackage{graphicx,subcaption,booktabs,tabularx,makecell,array}
\usepackage[ruled,vlined,linesnumbered]{algorithm2e}
\usepackage{newunicodechar}
\usepackage{enumitem}
\usepackage[hidelinks]{hyperref}
\usepackage{xurl}
\usepackage[labelfont=bf,labelsep=period,font=small]{caption}
\newunicodechar{−}{\ensuremath{-}}
\newunicodechar{±}{\ensuremath{\pm}}
\newunicodechar{≤}{\ensuremath{\le}}
\newunicodechar{≥}{\ensuremath{\ge}}
\newunicodechar{α}{\ensuremath{\alpha}}
\newunicodechar{✔}{\checkmark}
\newunicodechar{’}{'}
\SetAlgoCaptionSeparator{.}
\newcommand{\objref}[2]{\hyperref[#2]{\textbf{#1~\ref*{#2}}}}
\newenvironment{refs}{\begin{list}{}{\setlength{\leftmargin}{5mm}\setlength{\itemindent}{-5mm}\setlength{\itemsep}{2pt}\setlength{\parsep}{0pt}}\small}{\end{list}}
\title{\textbf{Intercity EV charging readiness in Türkiye: A countrywide scheduling and economic assessment}}
\author{Taner Cokyasar\\[2pt]
\small Texas A\&M University, Fermier Hall, 106 Ross St., College Station, TX, 77843, USA\\
\small \href{mailto:tcokyasar@tamu.edu}{tcokyasar@tamu.edu}, ORCID: \url{https://orcid.org/0000-0001-9687-6725}}
\date{}

\begin{document}
\maketitle

\begin{abstract}
This paper assesses whether Türkiye's 9,641 public charging stations licensed by the Energy Market Regulatory Authority (EPDK) can support intercity electric vehicle (EV) travel, and at what cost relative to internal combustion engine (ICE) vehicles. A charging scheduling framework based on a directed acyclic graph plans the charging stops of 3 million synthetic intercity trips on a 55,236-link highway network built from OpenStreetMap. Virtually all trips can be completed with the existing stations under both 10–100\% and 20–80\% state-of-charge windows; the 0.3\% feasibility that a 20–80\% window can produce arises only when trips must also start within the window, not from gaps in the network. Capacity, not coverage, is the binding constraint: with plug counts calibrated to EPDK socket statistics, uncoordinated charging builds long queues at the stations where trips converge, whereas coordinating charging times with reservations nearly eliminates waiting at today's adoption level of about 3\%. In the 350 kW scenario at 3\% adoption, the EV is 49.8\% cheaper than the ICE vehicle in operating costs (energy, refueling, detour, and charging time) and 7.4\% cheaper in total cost of ownership, an advantage that holds as long as the EV costs at most about 9\% more than a comparable ICE vehicle. Adding the cost of charging infrastructure makes the EV more expensive at low adoption with high-power chargers and cheaper as adoption grows.

\end{abstract}

\section{Introduction}\label{sec:1}

Republic of Türkiye stands at a pivotal moment in its transportation electrification trajectory. With 16.2 million registered passenger cars (TÜİK, \hyperlink{ref:tk2024}{2024}), which traveled 122.0 billion vehicle-kilometers on the motorways, state roads, and provincial roads in 2024 (KGM, \hyperlink{ref:kgm2025}{2025}), road transport accounts for a substantial share of the country's energy consumption and carbon emissions. Recognizing this, Türkiye has set ambitious targets for Electric Vehicle (EV) adoption, and EPDK-licensed operators run a rapidly growing network of public charging stations, of which the Şarj Gezgini platform listed 9,641 operated by 170 providers in July 2026 (Şarj Gezgini, \hyperlink{ref:arjgezgini2026}{2026}). Yet the fundamental question facing both policymakers and prospective EV buyers remains unanswered: can this emerging infrastructure support cost-competitive intercity travel, or does the detour and charging time overhead render long-distance EV trips economically impractical compared to conventional Internal Combustion Engine (ICE) vehicles?

Answering this question requires more than aggregate network statistics. It demands a vehicle-level analysis that routes individual trips through the existing charging network, determines where and when each vehicle should stop to recharge, resolves the inevitable conflicts that arise when multiple vehicles converge on the same station, and translates the resulting operational burdens into economic terms. This \emph{scheduling and assessment} perspective is fundamentally different from the \emph{siting} question of where to place new chargers~(Kuby \& Lim, \hyperlink{ref:kuby2005}{2005}; Li, Huang, \& Mason, \hyperlink{ref:li2016}{2016}) or the \emph{sizing} question of how many plugs each station should have~(Kınay, Gzara, \& Alumur, \hyperlink{ref:knay2023}{2023}; Upchurch, Kuby, \& Lim, \hyperlink{ref:upchurch2009}{2009}). While siting and sizing studies prescribe how infrastructure \emph{should} be deployed, a scheduling assessment reveals how well the infrastructure \emph{already deployed} serves the traveling public.

The existing literature on EV charging optimization, reviewed in Section~\ref{sec:2}, is dominated by studies set in North America and Western Europe, where mature highway networks and established charger registries provide rich input data. Emerging EV markets such as Türkiye, where charging networks are growing but unevenly distributed and where intercity distances can be substantial (e.g., over 1,400~km from Istanbul to Diyarbakır), have received far less attention. Moreover, most scheduling studies consider individual vehicles or small corridors rather than the millions of trips that define a country-level assessment.

This paper addresses that gap by applying the scalable charging schedule optimization framework of Cokyasar (\hyperlink{ref:cokyasar2026}{2026}), originally developed and validated for a U.S.~national-scale analysis, to the full Türkiye intercity highway network. The road network is constructed from OpenStreetMap motorway, trunk, and primary road data, yielding a connected graph of 31,286 nodes and 55,236 links. Intercity travel demand is represented by 3 million synthetic origin--destination trips distributed across 100,000 population-weighted locations spanning all 929 districts (\emph{ilçe}) of Türkiye. The charging infrastructure comprises 9,641 EPDK-licensed stations drawn from the Şarj Gezgini platform~(Şarj Gezgini, \hyperlink{ref:arjgezgini2026}{2026}). The methodological contribution of this paper lies not in the scheduling algorithm itself, which is detailed in~(Cokyasar, \hyperlink{ref:cokyasar2026}{2026}), but in the construction of Türkiye-specific inputs, the parametric adaptation to local conditions, and the policy-relevant interpretation of results for an emerging EV market.

The analysis yields several findings of policy relevance. Coverage is not the binding constraint: when a vehicle may start a trip above the upper bound of its State of Charge (SoC) window, virtually all of the 2.16 million trips that need charging are feasible under both 10--100\% and 20--80\% windows. Capacity is: with plug counts calibrated to EPDK statistics, uncoordinated charging produces long queues even at low adoption, whereas reservations keep waits below two minutes. EVs cost about half as much as ICE vehicles to operate, but their ownership advantage is narrow (7.4\% at 350 kW and 3\% adoption) and depends mainly on the EV-to-ICE price ratio.

The scheduling algorithm is that of Cokyasar (\hyperlink{ref:cokyasar2026}{2026}). The contributions of this paper are as follows:

\begin{enumerate}
\item
  To our knowledge, the first countrywide, vehicle-level assessment of intercity EV charging readiness for Türkiye, with documented inputs and trip counts at each processing step.
\item
  A congestion analysis with EPDK-calibrated plug counts and adoption levels of 3--100\%, comparing uncoordinated and reservation-based charging, together with a charging-gap analysis.
\item
  An ablation that measures how three modeling choices (restricting the starting SoC, discretizing SoC, and leaving charging time uncosted) change feasibility and charging behavior.
\item
  A three-tier cost comparison with discounted ownership costs, power-dependent infrastructure costs, break-even values of the price ratio and of the value of charging time, and sensitivity analyses.
\end{enumerate}

\section{Literature Review}\label{sec:2}

Prior work on EV charging spans infrastructure planning, vehicle-level routing, system-level scheduling, and economic evaluation. Rather than repeating a comprehensive survey, the reader is referred to Cokyasar (\hyperlink{ref:cokyasar2026}{2026}) for an extensive treatment, this section highlights the studies most relevant to Türkiye's context and identifies the gap that this paper fills. \objref{Table}{tab:1} positions representative studies along key dimensions.

\begin{table}[!htb]
\centering
\caption{Positioning of this study relative to representative prior work.}\label{tab:1}
\footnotesize
\begin{tabularx}{\textwidth}{>{\raggedright\arraybackslash}p{3.0cm} ccc >{\raggedright\arraybackslash}X >{\raggedright\arraybackslash}X >{\raggedright\arraybackslash}p{3.4cm}}
\toprule
& \multicolumn{3}{c}{\textbf{Problem}} & \multicolumn{2}{c}{\textbf{Modeling}} & \textbf{Approach} \\
\cmidrule(lr){2-4} \cmidrule(lr){5-6} \cmidrule(lr){7-7}
\textbf{Study} & \textbf{SI} & \textbf{ES} & \textbf{MC} & \textbf{Infrastructure} & \textbf{Capacity} & \textbf{Method / Scale} \\
\midrule
(Kuby \& Lim, \hyperlink{ref:kuby2005}{2005}) & \checkmark & & & Decision variable & None & Flow-refueling LP / Regional \\
(Upchurch et al., \hyperlink{ref:upchurch2009}{2009}) & \checkmark & & & Decision variable & Aggregate & Capacitated LP / Regional \\
(Li et al., \hyperlink{ref:li2016}{2016}) & \checkmark & & & Decision variable & None & Multi-period MIP / State \\
(Kınay et al., \hyperlink{ref:knay2023}{2023}) & \checkmark & & & Decision variable & Aggregate & Stochastic MIP / Regional \\
(Kullman et al., \hyperlink{ref:kullman2021}{2021}) & & \checkmark & & Fixed & None & Shortest path / Trip \\
(Bae \& Kwasinski, \hyperlink{ref:bae2012}{2012}) & & & \checkmark & Fixed & Analytical & M/M/c queuing / Station \\
(Lee et al., \hyperlink{ref:lee2019}{2019}) & & \checkmark & \checkmark & Fixed & Explicit & MIP / Facility \\
(Tang \& Zhang, \hyperlink{ref:tang2017}{2017}) & & & \checkmark & Fixed & Iterative & Dual decomp.\ / Facility \\
(Hagman et al., \hyperlink{ref:hagman2016}{2016}) & & & & Assumed & --- & TCO accounting / Fleet \\
(Borlaug et al., \hyperlink{ref:borlaug2020}{2020}) & & & & Assumed & --- & Levelized cost / Fleet \\
(Needell et al., \hyperlink{ref:needell2016}{2016}) & & & & Assumed & --- & Statistical / National \\
(Cokyasar, \hyperlink{ref:cokyasar2026}{2026}) & & \checkmark & \checkmark & Fixed (AFDC) & Iterative + FIFO & DAG DP / $\sim$2.7M vehicles, U.S. \\
\textbf{This paper} & & \checkmark & \checkmark & \textbf{Fixed (EPDK)} & \textbf{FIFO, reservation} & \textbf{DAG DP / $\sim$3M trips, Türkiye} \\
\bottomrule
\end{tabularx}

\smallskip
\raggedright\emph{SI: Siting/sizing, ES: En-route scheduling, MC: Multi-vehicle/capacity management. LP: linear program, MIP: mixed-integer program, DP: dynamic programming.}
\end{table}

\subsection{Infrastructure Planning versus Infrastructure Assessment}\label{sec:2_1}

The dominant strand of the EV charging literature treats infrastructure as a decision variable. Flow-refueling location models (Kuby \& Lim, \hyperlink{ref:kuby2005}{2005}; Upchurch et al., \hyperlink{ref:upchurch2009}{2009}) site stations to maximize intercepted traffic; coverage and budget formulations (Frade, Ribeiro, Gonçalves, \& Antunes, \hyperlink{ref:frade2011}{2011}; Li et al., \hyperlink{ref:li2016}{2016}) balance geographic reach against investment cost; and stochastic extensions (Kınay et al., \hyperlink{ref:knay2023}{2023}) hedge against demand uncertainty. All of these works answer the question \emph{where should we build?} In contrast, this paper asks \emph{what can we do with what has already been built?}, a diagnostic stance that is especially informative in markets where deployment decisions have already been made but their system-level adequacy has not been evaluated.

\subsection{Scheduling and Capacity at Scale}\label{sec:2_2}

At the single-vehicle level, en-route charging reduces to a resource-constrained shortest-path problem, solvable in seconds by dynamic programming (Kullman, Goodson, \& Mendoza, \hyperlink{ref:kullman2021}{2021}; Zhang, Chen, Yu, Zhu, \& Shi, \hyperlink{ref:zhang2018}{2018}). When multiple vehicles share finite plug supply, however, individually optimal plans can collectively violate station capacity. Queuing models (Bae \& Kwasinski, \hyperlink{ref:bae2012}{2012}), discrete-event simulation (Witt, \hyperlink{ref:witt2023}{2023}), and facility-level scheduling (Lee, Li, \& Low, \hyperlink{ref:lee2019}{2019}; Tang \& Zhang, \hyperlink{ref:tang2017}{2017}) each address aspects of this coupling but are limited in scale, typically to a single station, corridor, or fleet of tens to hundreds of vehicles. Cokyasar (\hyperlink{ref:cokyasar2026}{2026}), an arXiv preprint, proposed a three-stage framework (infeasibility pruning, per-vehicle DAG-based DP, and a congestion stage) that scales to millions of vehicles across a national network without requiring a commercial solver. The present paper applies that framework to Türkiye; the algorithmic contribution resides in Cokyasar (\hyperlink{ref:cokyasar2026}{2026}), and the contribution here lies in the data construction, parametric adaptation, and policy analysis for Türkiye.

\subsection{Economic Competitiveness of Long-Distance EV Travel}\label{sec:2_3}

Whether EVs can compete with ICE vehicles on long-distance trips is fundamentally an economic question. Fleet-level TCO studies~(Hagman, Ritzén, Stier, \& Susilo, \hyperlink{ref:hagman2016}{2016}) and levelized charger cost analyses~(Borlaug, Salisbury, Gerdes, \& Muratori, \hyperlink{ref:borlaug2020}{2020}) provide useful benchmarks but do not capture trip-level heterogeneity in detour burden, charging delay, and route geography. Statistical assessments based on travel surveys~(Chakraborty, Bunch, Lee, \& Tal, \hyperlink{ref:chakraborty2019}{2019}; Needell, McNerney, Chang, \& Trancik, \hyperlink{ref:needell2016}{2016}) estimate the share of vehicle-miles serviceable by EVs but do not optimize charging decisions. The framework used here evaluates costs at the individual path level, decomposing the EV--ICE differential into energy, detour time, waiting time, vehicle ownership, and infrastructure amortization with detour and wait costs produced endogenously by the scheduler rather than assumed exogenously.

\subsection{EV Readiness in Emerging Markets}\label{sec:2_4}

Large-scale charging analyses have concentrated on North America and Western Europe, where mature registries (e.g., the U.S. Alternative Fuels Data Center, the European Alternative Fuels Observatory) and household travel surveys provide input data (IEA, \hyperlink{ref:iea2024}{2024}). For Türkiye, Gulbahar, Sutcu, Almomany, and Ibrahim (\hyperlink{ref:gulbahar2023}{2023}) forecast the EV fleet and locate charging stations on state highways to meet intercity travel demand, and Gönül, Duman, and Güler (\hyperlink{ref:gnl2024}{2024}) propose a weighted-sum framework for siting stations along highways and demonstrate it on Turkish motorways. Studies of adoption find that vehicle taxes are the main factor in consumers' decision to switch to EVs and that more charging stations raise willingness to switch (Şen, Yiğiter, \& Özcan, \hyperlink{ref:en2023}{2023}), and Durmuş Şenyapar, Akıl, and Dokur (\hyperlink{ref:durmuenyapar2023}{2023}) review the Turkish research on EV purchase intention; EPDK projects the growth of the EV fleet and of charging infrastructure to 2040 (EPDK, \hyperlink{ref:epdk2026b}{2026b}). These studies plan new stations or study adoption; they do not simulate individual intercity trips on the existing station network. Türkiye's network, with 9,641 public stations listed in July 2026 (Şarj Gezgini, \hyperlink{ref:arjgezgini2026}{2026}), has evolved around metropolitan corridors, leaving open the question of whether intercity coverage is sufficient. To our knowledge, no prior work has conducted a countrywide, vehicle-level assessment of intercity EV charging readiness for Türkiye.

\section{Methodology}\label{sec:3}

This section is organized into three subsections. We provide an overview of the scheduling framework in Section~\ref{sec:3_1}. In Section~\ref{sec:3_2}, we then review how the Türkiye-specific data are constructed. Finally, we define the framework of economic evaluation conducted in Section~\ref{sec:3_3}.

\subsection{Scheduling Framework Overview}\label{sec:3_1}

We employ the three-stage EV charging scheduling framework developed by Cokyasar (\hyperlink{ref:cokyasar2026}{2026}), an arXiv preprint; the reader is referred to that work for the full mathematical formulation. Here we provide a self-contained summary of the elements necessary to understand the Türkiye application, and \objref{Algorithm}{alg:1} gives the pseudocode of the per-vehicle step. Three features of the per-vehicle step, described below, matter for the results: a vehicle's starting SoC may exceed the charging target, SoC is tracked exactly rather than in 1\% buckets, and the time spent charging is part of the cost.

The problem can be described as follows. Consider a fleet of EVs, each assigned a predetermined intercity route through the road network. A route consists of an ordered sequence of directed links, obtained from a static shortest-path computation between the vehicle's origin and destination. Along each route, a set of charging stations with known locations and finite plug counts is available. After traversing a link, a vehicle may choose to detour from the link's endpoint to one of the \(K\) nearest charging stations, recharge its battery to a target State of Charge (SoC), and return to the same endpoint before continuing. The central scheduling question is: for each vehicle, after which links should it detour to charge? The answer must respect two families of constraints. First, the vehicle's SoC must remain within an operational window \(\left\lbrack \underline{S},\overline{S} \right\rbrack\) at every point along the route, where \(\underline{S}\) is the minimum acceptable SoC (a safety buffer against deep discharge) and \(\overline{S}\) is the maximum to which the battery is charged at each stop. A vehicle starts its trip with SoC \(S_{v}^{0}\), which may exceed \(\overline{S}\): a vehicle charged at home to 100\% is not treated as if it had been charged only to \(\overline{S}\). Second, the number of vehicles simultaneously occupying plugs at any station must not exceed the station's capacity. The objective is to minimize the total cost incurred by all vehicles, comprising the time cost of round-trip detours to and from charging stations, the time cost of charging, and the waiting cost when a vehicle arrives at a station whose plugs are all occupied.

This joint scheduling problem can be formulated as a mixed-integer program (MIP), but the resulting model is computationally intractable at the scale of the Türkiye application: with over two million vehicles, routes averaging hundreds of links each, and nearly ten thousand stations, the number of binary decision variables exceeds \(10^{8}\) before capacity-coupling variables are introduced. However, the problem has a key structural property: if station capacity constraints are temporarily ignored, each vehicle's scheduling subproblem becomes independent and can be solved optimally in isolation. The capacity constraint is the sole source of inter-vehicle coupling. This decomposable structure motivates a three-stage solution pipeline described next.

\emph{Stage~1: Infeasibility Pruning.} Before attempting to optimize charging schedules, the framework identifies vehicles that cannot possibly complete their journey without violating the minimum SoC constraint, regardless of how charging stops are arranged. This occurs when a vehicle's route contains a stretch of consecutive links between two charging opportunities whose cumulative energy consumption exceeds the usable SoC window \(\overline{S} - \underline{S}\). A forward-pass reachability heuristic performs this check efficiently by simulating a greedy charging policy along each route: at every link with charger access, it evaluates whether the vehicle has enough remaining SoC to reach the next charger; if not, it charges immediately; if charging is impossible (the charger is unreachable even from the current SoC), the vehicle is declared infeasible and excluded. This pruning step avoids wasting computational resources on vehicles for which no feasible schedule exists.

\emph{Stage 2: DAG Shortest-Path Solver.} For each surviving vehicle, the framework constructs a Directed Acyclic Graph (DAG) whose nodes are the vehicle's charging opportunities, plus the origin and the destination. An opportunity \((i,k)\) combines the end of link \(i\) of the route with one of its candidate stations \(k\): the nearest station and up to \(K - 1\) further stations within a distance \(\Delta\) of the link end. An arc from an opportunity \((h,m)\) to a later opportunity \((i,k)\) exists if the vehicle, after charging to \(\overline{S}\) at \((h,m)\), can reach \((i,k)\) without falling below \(\underline{S}\); arcs from the origin use \(S_{v}^{0}\) instead of \(\overline{S}\). Because every stop ends at the same target \(\overline{S}\), the SoC on arrival at each opportunity is known exactly, and no discretization of SoC is needed. The cost of a stop is the time cost of the detour plus the time cost of charging from the arrival SoC to \(\overline{S}\) at constant power. Because the opportunities follow the route order, the cheapest schedule is found by a single forward Dynamic Programming (DP) pass (\objref{Algorithm}{alg:1}). A common alternative, a DAG over (link, SoC bucket) pairs with 1\% buckets, is less accurate: rounding each link's energy use to whole buckets overstates the energy of the many short links of a road network and inflates the number of stops. All vehicles are solved independently in parallel across available processor cores.

\begin{algorithm}[!htb]
\caption{Capacity-free charging schedule of one vehicle (charging-opportunity DAG)}\label{alg:1}
\small
\KwIn{route links $1,\ldots,n$ with energy use $e_l$ (\% SoC); candidate stations $k$ at the end of each link $i$ with detour time $t_{ik}$ (hours) and detour energy $d_{ik}$ (\% SoC); starting SoC $S_v^0$; operating bounds $\underline{S}$, $\overline{S}$; charger power $W$ (kW); battery capacity $Q_v$ (kWh); time cost rates $c^{\tau}$ and $c^{\chi}$ (TL/h)}
\KwOut{the cheapest set of charging stops}
\textbf{Notation:} $U(h,i)=\sum_{l=h+1}^{i} e_l$ is the SoC used between the ends of links $h$ and $i$ (link 0 ends at the origin); $s_{hm}$ is the SoC on arrival at opportunity $(h,m)$ from its predecessor; $\text{cost}(i,k)=\infty$ until set; $\text{stop}(h,m)=c^{\tau}t_{hm}+c^{\chi}Q_v(\overline{S}-s_{hm})/(100W)$\;
\ForEach{opportunity $(i,k)$ in route order}{
  \lIf{$S_v^0-U(0,i)-d_{ik}/2\ge\underline{S}$}{$\text{cost}(i,k)\leftarrow 0$, $\text{pred}(i,k)\leftarrow$ origin}
  \ForEach{earlier opportunity $(h,m)$ with $\overline{S}-d_{hm}/2-U(h,i)-d_{ik}/2\ge\underline{S}$}{
    \If{$\text{cost}(h,m)+\text{stop}(h,m)<\text{cost}(i,k)$}{
      $\text{cost}(i,k)\leftarrow\text{cost}(h,m)+\text{stop}(h,m)$, $\text{pred}(i,k)\leftarrow(h,m)$\;
    }
  }
}
$(i',k')\leftarrow$ the opportunity with the lowest $\text{cost}(i,k)+\text{stop}(i,k)$ from which the destination is reached at or above $\underline{S}$\;
\Return the stops on the predecessor chain from $(i',k')$ back to the origin\;
\end{algorithm}

\emph{Stage 3: Congestion Resolution.} The DAG solver produces individually optimal schedules, but when these are combined across all vehicles, some stations may be assigned more simultaneous users than they have plugs. In the baseline (no coordination), each vehicle keeps its capacity-free schedule, and a First-In, First-Out (FIFO) event simulation processes all charging sessions in time order: a vehicle that finds all plugs busy waits, and its delay is carried to its later stops. This gives every vehicle a capacity-feasible schedule with an explicit waiting cost. As a coordinated alternative, sequential reservation plans vehicles one at a time in order of departure; each vehicle solves a time-dependent version of \objref{Algorithm}{alg:1} against the plug reservations of earlier vehicles, so that waiting is part of its cost, and then reserves its own sessions. Reservation represents an idealized booking service that all drivers use and honor. Cokyasar (\hyperlink{ref:cokyasar2026}{2026}) reports the accuracy of both methods against an exact mixed-integer program on small instances; an iterative station-penalty heuristic, also tested there, performed worse than no coordination in most cases and is not used.

\subsection{Türkiye-Specific Data Construction}\label{sec:3_2}

Applying the scheduling framework to Türkiye required constructing three datasets from scratch: an intercity road network, a synthetic travel demand, and a georeferenced charging station registry. This subsection details each.

The intercity highway network was extracted from OpenStreetMap (OSM) (OpenStreetMap contributors, \hyperlink{ref:openstreetma2026}{2026}) via the Geofabrik regional extract for Türkiye (Geofabrik, \hyperlink{ref:geofabrik2026}{2026}). Road segments classified as motorway, trunk, or primary (including their respective \_link variants) were retained; lower-class roads were excluded to restrict the network to intercity-relevant facilities. Segment endpoints within 100 m of each other were merged using a KD-tree spatial snap to create intersection nodes. The largest connected component was retained, yielding 31,286 nodes and 55,236 directed links. Free-flow speeds were assigned from OSM maxspeed tags where available, with highway-class defaults (motorway: 120 km/h, trunk: 100 km/h, primary: 80 km/h) applied otherwise. All coordinates were projected to EPSG:32636 (UTM Zone 36N). The network was stored in the POLARIS supply database format (Auld et al., \hyperlink{ref:auld2016}{2016}), and all routes were computed with the static shortest-path router of AequilibraE (AequilibraE developers, \hyperlink{ref:aequilibraed2026}{2026}).

No publicly available intercity origin--destination trip table exists for Türkiye. We constructed a synthetic demand of 3 million personal vehicle trips representing three days of intercity travel, that is, one million trips per day. KGM does not publish trip counts or an origin--destination table; it publishes vehicle-kilometers by road class and vehicle type. The one million trips per day is therefore our assumption, which we checked against KGM's statistics in vehicle-kilometers. The trips that need charging alone amount to about 228 billion vehicle-kilometers per year, 1.9 times the 122.0 billion car vehicle-kilometers on KGM's network in 2024 (KGM, \hyperlink{ref:kgm2025}{2025}, Table UE04), because drawing origins and destinations by population across the country makes them long (866 km on average). The synthetic demand is therefore a stress test weighted toward long trips: feasibility and congestion results are conservative, whereas per-trip costs do not depend on total volume. Origin and destination locations were drawn from 100,000 population-weighted points distributed across all 929 administrative districts (\emph{ilçe}), with district populations from TÜİK (TÜİK, \hyperlink{ref:tk2025a}{2025a}) allocated using GADM Level 2 boundaries (GADM, \hyperlink{ref:gadm2022}{2022}). Trip start times were distributed across three days using an hourly profile that peaks during morning and evening periods. Origins were sampled proportionally to district population; destinations were drawn from a different district to ensure intercity character. All trips were routed through the network, producing link-level trajectories.

Not every synthetic trip enters the scheduler. Of the 3,000,000 trips, 187,046 (6.2\%) could not be routed because their origin or destination lies on a link that the rest of the directed network cannot reach or leave; they are not predominantly short (median straight-line distance of 383 km, against 448 km for routed trips), so their removal does not favor easy trips. Of the 2,812,954 routed trips, 651,704 need no charging. Of the remaining 2,161,250, 3 cannot be completed with any charging plan and are removed in Stage 1, which leaves 2,161,247 for the scheduler.

Station locations were obtained from Şarj Gezgini (Şarj Gezgini, \hyperlink{ref:arjgezgini2026}{2026}), a third-party platform that lists the public charging stations of operators licensed by EPDK; it is not an EPDK database. The data, retrieved on July 13, 2026, contain 9,641 stations operated by 170 networks (largest: ZES with 1,827 stations, Eşarj with 745, Trugo with 649). Geographic coordinates (WGS 84) were projected to EPSG:32636 for consistency with the road network. Because the platform does not report plug counts per station, the case study brackets the unknown capacity with 50 and 5 plugs per station. EPDK, however, reports national socket totals: in June 2026, there were 19,663 DC and 25,434 AC public sockets (EPDK, \hyperlink{ref:epdk2026a}{2026a}), about two DC sockets per listed station. Section~\ref{sec:4_4} uses this calibrated capacity to study realistic adoption levels. Station-to-link proximity was computed using a BallTree index with Manhattan distance in projected coordinates.

\subsection{Economic Evaluation Framework}\label{sec:3_3}

Following Cokyasar (\hyperlink{ref:cokyasar2026}{2026}), economic outcomes are evaluated through a three-tier analysis computed at the individual path level and then aggregated. Each tier progressively expands the cost boundary from operational expenses to full system costs. All distances are in kilometers and all costs in Turkish lira (TL). The notation used in the tier equations is as follows. Let \(\mathcal{V}\) denote the set of vehicles, one for each trip that enters the scheduler (Section~\ref{sec:3_2}). For each vehicle \(v \in \mathcal{V}\): \(D_{v}\) is the original trip distance (km), \(D_{v}^{\delta}\) is the total round-trip detour distance (km) incurred by charging stops, \(C_{v}^{\tau}\) is the total detour time cost (TL), valued at the rate \(c^{\tau}\), \(C_{v}^{\chi}\) is the total charging time cost (TL), \(C_{v}^{\omega}\) is the total waiting cost (TL), valued at the rate \(c^{\omega}\), and \(\rho_{v}\) is the vehicle's discharge rate (Wh/km). \(C_{v}^{\tau}\), \(C_{v}^{\chi}\), and \(C_{v}^{\omega}\) are outputs of the scheduling optimization, not assumed parameters. The remaining symbols are shared parameters: \(P^{\gamma}\) is the gasoline price (TL/liter), \(P^{\xi}\) is the electricity price (TL/kWh), which differs for AC and DC charging (\objref{Table}{tab:2}), \(\eta^{\gamma}\) is ICE fuel economy (km/liter), and \(R^{\gamma}\) and \(\theta^{\gamma}\) are the ICE driving range between refueling stops (km) and the duration of a refueling stop (h), valued at the same rate \(c^{\chi}\) as EV charging time.

Tier 1 (Operational) compares ICE fuel and refueling time costs against EV electricity, detour, charging time, and waiting costs:

\begin{align}
C_v^{\gamma,1} &= \frac{D_v P^{\gamma}}{\eta^{\gamma}} + c^{\chi}\theta^{\gamma}\left\lfloor \frac{D_v}{R^{\gamma}} \right\rfloor \label{eq:1}\\
C_v^{\xi,1} &= \frac{P^{\xi}\left(D_v + D_v^{\delta}\right)\rho_v}{1000} + C_v^{\tau} + C_v^{\chi} + C_v^{\omega} \label{eq:2}
\end{align}

Tier 2 (Total Cost of Ownership) adds vehicle ownership and maintenance costs, both charged to a trip in proportion to its distance. Let \(V^{\gamma}\) denote the ICE vehicle purchase price, \(\alpha\) the EV-to-ICE price ratio, \(M^{\gamma}\) and \(M^{\xi}\) the per-km maintenance costs for ICE and EV respectively, \(L\) the economic life of the vehicle in years, \(r\) the real discount rate, \(\phi^{\gamma}\) and \(\phi^{\xi}\) the residual values after \(L\) years as shares of the purchase price, \(\beta\) the cost of an EV battery replacement in year \(L^{b}\) (zero in the baseline), and \(A\) the annual distance driven by a car (km/yr). The ownership cost per km, \(a^{\gamma}\) and \(a^{\xi}\), is the purchase price net of the discounted residual value (plus the discounted battery replacement for the EV), annualized with the capital recovery factor \(\mathrm{CRF}(r,L) = r(1 + r)^{L}/\left\lbrack (1 + r)^{L} - 1 \right\rbrack\) and divided by \(A\):

\begin{align}
C_v^{\gamma,2} &= C_v^{\gamma,1} + \left(a^{\gamma} + M^{\gamma}\right)D_v, \quad a^{\gamma} = \frac{V^{\gamma}\left[1 - \phi^{\gamma}(1+r)^{-L}\right]\mathrm{CRF}(r,L)}{A} \label{eq:3}\\
C_v^{\xi,2} &= C_v^{\xi,1} + \left(a^{\xi} + M^{\xi}\right)\left(D_v + D_v^{\delta}\right), \quad a^{\xi} = \frac{\left[\alpha V^{\gamma}\left(1 - \phi^{\xi}(1+r)^{-L}\right) + \beta(1+r)^{-L^{b}}\right]\mathrm{CRF}(r,L)}{A} \label{eq:4}
\end{align}

Tier 3 (Infrastructure) adds amortized charging infrastructure capital and operating costs to the EV side. Let \(\mathcal{P}\) denote the set of charger power levels in the network and \(Z_{p}\) the number of plugs of power level \(p \in \mathcal{P}\), \(\kappa_{p}\) the installed capital cost per plug at that power, \(L_{c}\) the charger useful life in years, \(\mu\) the annual operations and maintenance rate as a fraction of capital cost, \(E\) the annual energy delivered by the network (kWh/yr), and \(E_{v}\) the energy vehicle \(v\) charges en route (kWh). In each scenario, all public plugs have the scenario's charger power \(W\), so \(\mathcal{P} = \{ W\}\) and \(Z_{W}\) is the number of AC plugs reported by EPDK if \(W \leq 22\) kW and of DC plugs otherwise. The annualized cost of the network is allocated per kWh delivered:

\begin{equation}
C_v^{\xi,3} = C_v^{\xi,2} + \frac{\sum_{p\in\mathcal{P}} Z_p \kappa_p \left[\mathrm{CRF}(r,L_c) + \mu\right]}{\mathrm{E}}\, E_v \label{eq:5}
\end{equation}

The ICE cost at Tier 3 is unchanged (\(C_{v}^{\gamma,3} = C_{v}^{\gamma,2}\)). Because retail charging prices already recover part of the infrastructure cost, Tier 3 is an upper bound on the cost of en-route charging. At each tier \(g \in \{ 1,2,3\}\), we report savings as a percentage of the total ICE cost over all trips, \(100 \times \left( \sum_{v \in \mathcal{V}}^{}C_{v}^{\gamma,g} - \sum_{v \in \mathcal{V}}^{}C_{v}^{\xi,g} \right)/\sum_{v \in \mathcal{V}}^{}C_{v}^{\gamma,g}\), where positive values indicate that EVs are cheaper. Under uncoordinated charging, waits measure demand that stations cannot serve rather than a cost a traveler would accept (Section~\ref{sec:4_4}), so the economic results exclude the waiting cost \(C_{v}^{\omega}\), which is reported separately; under reservation-based coordination, including it changes Tier 1 by less than 0.1 percentage points.

\section{Numerical Experiments}\label{sec:4}

\subsection{Experimental Design}\label{sec:4_1}

The travel demand comprises 3 million synthetic intercity personal vehicle trips representing three days of travel across Türkiye, as described in Section~\ref{sec:3_2}. All trips are routed through the intercity highway network (31,286 nodes, 55,236 links) using AequilibraE's static shortest-path algorithm, producing link-level trajectories with per-link distances, travel times, and sequences. EV-specific characteristics (range, discharge rate, and starting SoC) are assigned stochastically to each trajectory: range \(R_{v}\) is drawn from a uniform distribution over \(\left\lbrack \underline{R},\overline{R} \right\rbrack\), discharge rate \(\rho_{v}\) from a triangular distribution on \(\left\lbrack \underline{\rho},\overline{\rho} \right\rbrack\) with mode \(\underline{\rho}\), and starting SoC \(S_{v}^{0}\) from a triangular distribution on \(\left\lbrack {\underline{S}}^{0},{\overline{S}}^{0} \right\rbrack\) with mode \({\overline{S}}^{0}\). The battery capacity is \(Q_{v} = R_{v}\rho_{v}/1000\) (kWh).

The charging infrastructure consists of the 9,641 EPDK-licensed stations from the Şarj Gezgini platform (Şarj Gezgini, \hyperlink{ref:arjgezgini2026}{2026}). The effect of charger technology is examined by varying charger power \(W\) across four levels (5, 22, 50, and 350 kW) representing slow AC, fast AC, standard DC fast charging, and ultra-fast DC charging, respectively. Because plug counts per station are not published, the case study (Section~\ref{sec:4_2}) brackets capacity with 50 and 5 plugs per station, and the adoption analysis (Section~\ref{sec:4_4}) uses the 19,663 DC plugs reported by EPDK, about two per station. The detour speed \(\psi\) governs the time cost of the Manhattan distance detour from a link endpoint to a candidate station.

\objref{Table}{tab:2} summarizes all parameter values, organized by tier, with their units, baseline values, and references. All monetary values are denominated in Turkish Lira (TL), distances in kilometers, and energy consumption in Wh/km, reflecting Türkiye-specific conditions.

\begin{table}[!htb]
\centering
\caption{Parameter values used in the experiments.}\label{tab:2}
\footnotesize
\begin{tabularx}{\textwidth}{l >{\raggedright\arraybackslash}X >{\raggedright\arraybackslash}p{3.3cm} >{\raggedright\arraybackslash}p{3.2cm}}
\toprule
\textbf{Symbol} & \textbf{Meaning (unit)} & \textbf{Value / Range} & \textbf{Reference} \\
\midrule
\multicolumn{4}{l}{\emph{Simulation parameters}} \\
$\underline{S}^0$, $\overline{S}^0$ & Starting SoC bounds (\%) & 80 -- 100 & Assumed \\
$c^{\tau}$ & Detour cost rate (TL/h) & 100 & (Numbeo, \hyperlink{ref:numbeo2025}{2025}) \\
$c^{\omega}$ & Waiting cost rate (TL/h) & 70 & (Numbeo, \hyperlink{ref:numbeo2025}{2025}) \\
$K$ & Candidate stations per link endpoint & 3 & -- \\
$\underline{R}$, $\overline{R}$ & EV range bounds (km) & 250 -- 550 & (EV Database, \hyperlink{ref:evdatabase2025}{2025}) \\
$\underline{S}$, $\overline{S}$ & Operational SoC (\%) bounds & 10 -- 100 & (Keil \& Jossen, \hyperlink{ref:keil2017}{2017}) \\
$W$ & Charger power (kW) & 5, 22, 50, 350 & (ZES, \hyperlink{ref:zes2026}{2026}) \\
$c^{\chi}$ & Charging time cost rate (TL/h) & 70 & (Numbeo, \hyperlink{ref:numbeo2025}{2025}) \\
$\Delta$ & Maximum distance from the link end of the 2nd to $K$th candidate station (km) & 25 & Assumed \\
$\psi$ & Assumed detour speed (km/h) & 90 & Assumed \\
$\underline{\rho}$, $\overline{\rho}$ & Discharge rate bounds (Wh/km) & 150 -- 250 & (EV Database, \hyperlink{ref:evdatabase2025}{2025}) \\
\multicolumn{4}{l}{\emph{Tier~1 economic parameters (operational)}} \\
$P^{\gamma}$ & Gasoline price (TL/liter) & 75.0 & (Petrol Ofisi, \hyperlink{ref:petrolofisi2026}{2026}) \\
$P^{\xi}$ & Electricity price, AC ($W\le 22$ kW) and DC (TL/kWh) & 10.0, 14.5 & (Trugo, \hyperlink{ref:trugo2026}{2026}; ZES, \hyperlink{ref:zes2026}{2026}) \\
$\eta^{\gamma}$ & ICE fuel economy (km/liter) & 14.0 & (IEA, \hyperlink{ref:iea2024}{2024}) \\
$R^{\gamma}$ & ICE range between refueling stops (km) & 700 & Assumed \\
$\theta^{\gamma}$ & ICE refueling stop duration (h) & 0.167 (10 min) & Assumed \\
\multicolumn{4}{l}{\emph{Tier~2 economic parameters (total cost of ownership)}} \\
$A$ & Annual distance driven by a car (km/yr) & 13,859 & (TÜİK, \hyperlink{ref:tk2025b}{2025b}) \\
$L$ & Vehicle economic life (years) & 10 & Assumed \\
$M^{\gamma}$ & ICE maintenance cost (TL/km) & 0.35 & (Numbeo, \hyperlink{ref:numbeo2025}{2025}) \\
$M^{\xi}$ & EV maintenance cost (TL/km) & 0.12 & (Numbeo, \hyperlink{ref:numbeo2025}{2025}) \\
$V^{\gamma}$ & ICE vehicle purchase price (TL) & 2,545,500 & (Toyota Türkiye, \hyperlink{ref:toyotatrkiye2026}{2026}) \\
$\alpha$ & EV purchase price ratio (EV/ICE) & 1.0 & (Togg, \hyperlink{ref:togg2026}{2026}; Toyota Türkiye, \hyperlink{ref:toyotatrkiye2026}{2026}) \\
$r$ & Real discount rate (fraction) & 0.08 & Assumed \\
$\phi^{\gamma}$, $\phi^{\xi}$ & Residual value after $L$ years, ICE and EV (share of price) & 0.30, 0.25 & Assumed \\
$\beta$ & EV battery replacement cost, incurred in year $L^{b}$ (TL) & 0 & Assumed \\
\multicolumn{4}{l}{\emph{Tier~3 economic parameters (infrastructure)}} \\
$L_c$ & Charger useful life (years) & 12 & Assumed \\
$\kappa_p$ & Installed cost per plug: 5--22 kW, 50 kW, 350 kW (TL) & 270,000; 3,000,000; 9,000,000 & (Nicholas, \hyperlink{ref:nicholas2019}{2019}), at 47 TL/USD \\
$Z_p$ & Public plugs (sockets) in operation, AC and DC & 25,434; 19,663 & (EPDK, \hyperlink{ref:epdk2026a}{2026a}) \\
$\mu$ & Charger annual O\&M (fraction of capital cost) & 0.08 & Assumed \\
\bottomrule
\end{tabularx}
\end{table}

\objref{Figure}{fig:1} shows the geographic distribution of the 9,641 charging stations alongside the intercity highway network. Station density is highest along the Istanbul--Ankara--Izmir triangle and the Mediterranean coastal corridor, with sparser coverage in the eastern provinces. This uneven spatial distribution, which reflects the market-driven rollout of charging infrastructure around major metropolitan demand centers, determines where gaps in coverage remain.

\begin{figure}[!htb]
\centering
\includegraphics[width=\textwidth]{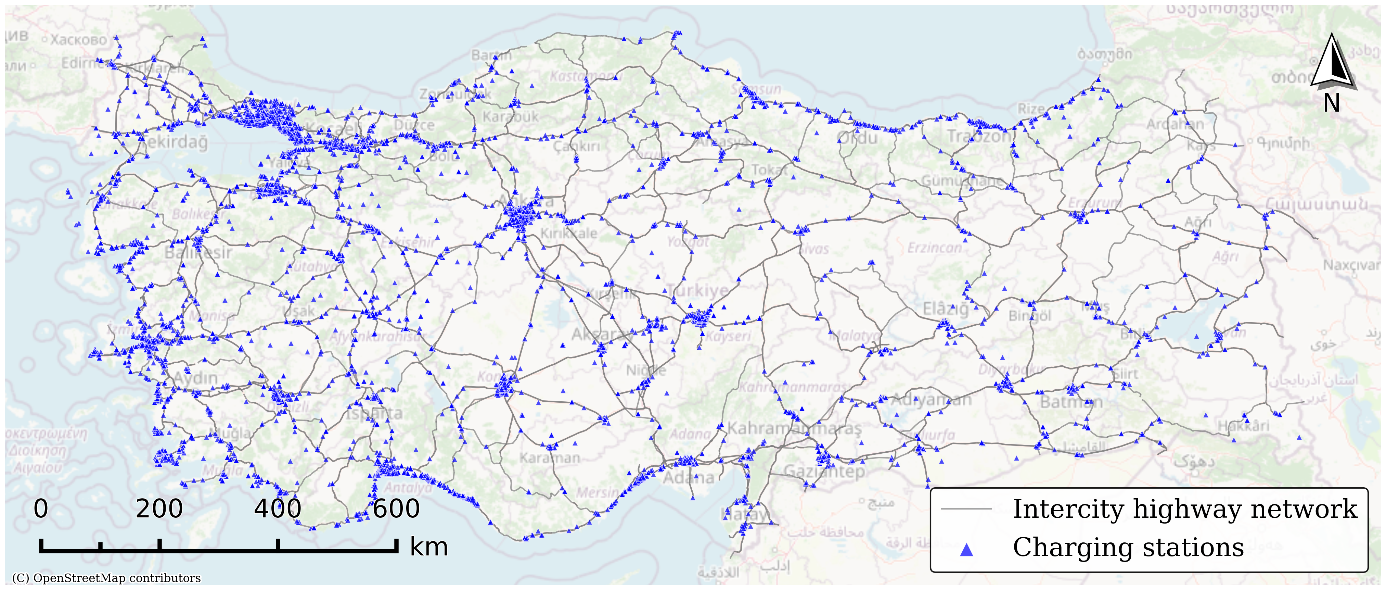}
\caption{Geographic distribution of 9,641 EPDK-licensed charging stations (blue triangles) and the intercity highway network (55,236 links, gray lines) across Türkiye.}\label{fig:1}
\end{figure}

Measured as the Manhattan distance from the end of each link to the nearest station, 42.8\% of the network's 71,925 km lies more than 10 km from a station, 13.8\% more than 25 km, and 1.3\% more than 50 km. The longer gaps are concentrated in eastern and northeastern Anatolia (\objref{Figure}{fig:2}); a vehicle can bridge them with the energy it carries, so they lengthen detours rather than prevent trips (Section~\ref{sec:4_3}).

\begin{figure}[!htb]
\centering
\includegraphics[width=\textwidth]{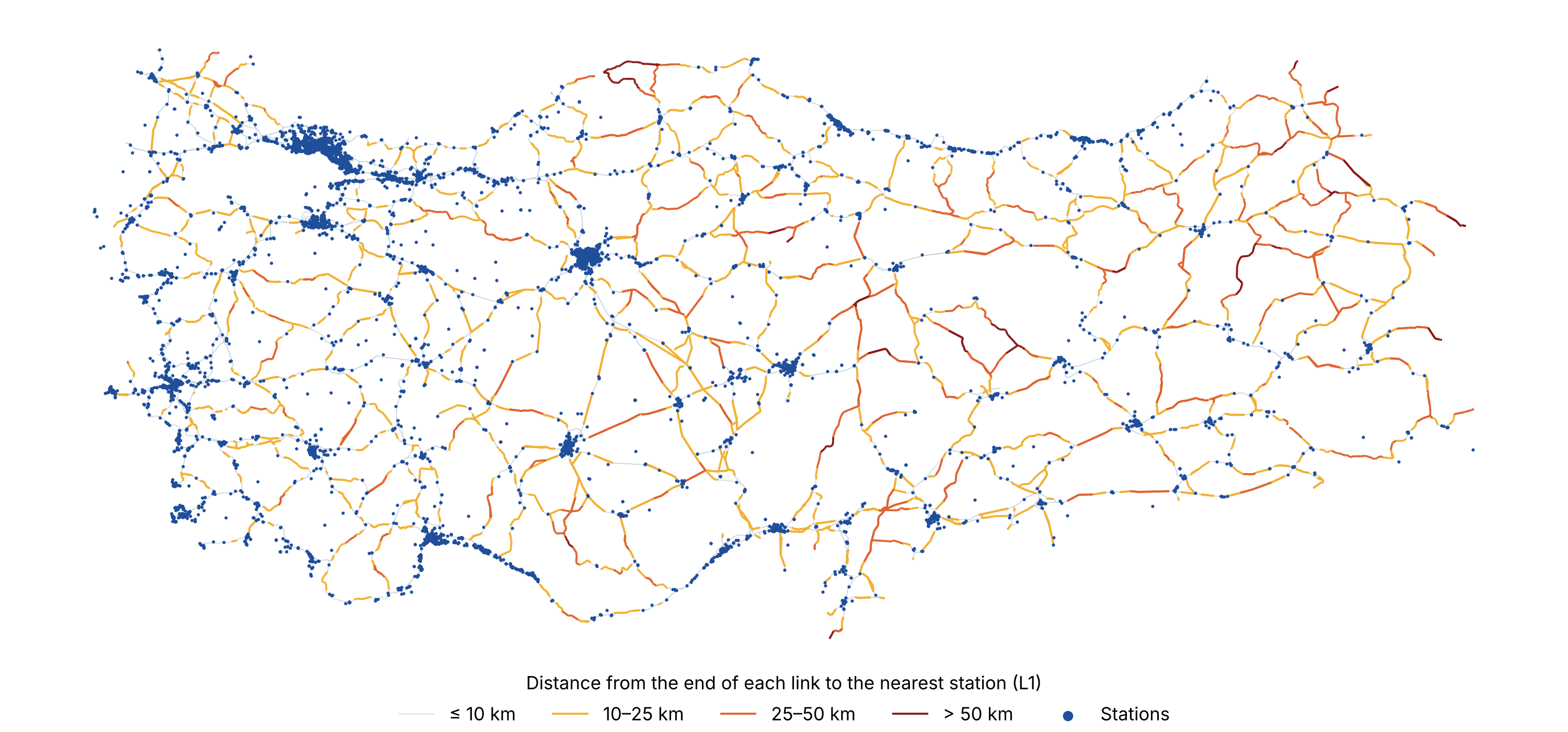}
\caption{Distance from the end of each link of the intercity network to the nearest charging station (Manhattan distance).}\label{fig:2}
\end{figure}

\subsection{Case Study}\label{sec:4_2}

The case study crosses four charger power levels (5, 22, 50, and 350 kW) with two plug counts that bracket the unknown station capacity: 50 plugs per station (482,050 in total) and 5 plugs per station (48,205), both above the roughly two DC plugs per station reported by EPDK (Section~\ref{sec:4_4}). All eight scenarios use the 9,641 stations, the 3 million synthetic trips, full electrification, the 10--100\% SoC window, and uncoordinated FIFO charging. \objref{Table}{tab:3} reports the results and the computing time of each run.

Virtually all trips are feasible: of the 2,161,247 trips that enter the scheduler, at most 7 cannot be completed. Because schedules are planned without regard to capacity, the plug count affects only waiting. A trip makes about 2.4--2.5 stops, detours add 0.26--0.59\% to its distance (1.5--3.4 minutes), and a stop lasts from 7.2 minutes at 350 kW to 497 minutes at 5 kW.

Waiting is where capacity binds: under uncoordinated charging, mean waits per trip range from 4.5 days (350 kW, 50 plugs) to several years (5 kW, 5 plugs). These values indicate demand that stations cannot serve within the simulated period rather than a cost a traveler would accept; Section~\ref{sec:4_4} examines realistic adoption and coordination.

EVs are cheaper to operate in every scenario, with Tier 1 savings of 35.2\% at 5 kW, 58.7\% at 22 kW, 47.1\% at 50 kW, and 49.8\% at 350 kW. AC charging (22 kW or less) costs 10 TL per kWh against 14.5 TL for DC, and the time cost of charging grows as power falls. Gasoline at 75 TL per liter costs 5.4 TL per km, against about 2.9 TL per km for DC electricity at 200 Wh per km.

The break-even value of charging time, at which the EV and ICE trips cost the same, is about 150 TL per hour at 5 kW, 680 at 22 kW, 1,220 at 50 kW, and 13,100 at 350 kW at Tier 1 (130, 580, 990, and 10,800 at Tier 2). Against the assumed 70 TL per hour, the EV remains cheaper for any plausible valuation of time at 22 kW and above, but not at 5 kW for a traveler who values time at about twice the assumed rate.

Ownership costs narrow the advantage to 4.6--9.0\% at Tier 2 because vehicles are expensive relative to fuel in Türkiye: annualized over a 10-year life at an 8\% discount rate, net of residual value, owning an ICE vehicle costs 23.6 TL per km, 4.4 times its fuel cost. Even with the straight-line amortization of the U.S. analysis of Cokyasar (\hyperlink{ref:cokyasar2026}{2026}), this ratio is 2.8 in Türkiye against about 1.5 in the United States. Tier 2 therefore depends mainly on the EV-to-ICE price ratio α. At November 2026 list prices, the entry Togg T10X costs 1.91 million TL and a gasoline Toyota Corolla 2.55 million TL (Togg, \hyperlink{ref:togg2026}{2026}; Toyota Türkiye, \hyperlink{ref:toyotatrkiye2026}{2026}), largely because of the lower special consumption tax on EVs; we conservatively set α = 1.0. The EV keeps its advantage up to a break-even α of 1.06--1.11 (1.09 in the reference scenario of 350 kW and 3\% adoption with reservation; \objref{Table}{tab:6}); across α from 0.75 to 1.4, Tier 2 savings in that scenario range from 28.2\% to −25.8\%. Varying each economic input by ±20\% in the reference scenario confirms this ranking: Tier 2 savings range from 24.0\% to −9.2\% with α and move by 6.8, 3.6, 3.5, and 3.4 percentage points with the gasoline price, electricity price, ICE price, and annual distance, while the value of charging time has no visible effect at 350 kW.

Infrastructure decides the full-system comparison at low adoption. With plug costs that rise from 0.27 million TL (22 kW AC) to 3.0 million TL (50 kW) and 9.0 million TL (350 kW), spread over each scenario's charging energy, Tier 3 savings are positive at full electrification (4.6--9.0\%; \objref{Table}{tab:3}). At 3\% adoption with reservation, the EV and ICE vehicle cost the same at 50 kW, whereas a 350 kW network makes the EV 11.4\% more expensive (\objref{Table}{tab:6}; 0.6\% and −11.3\% under FIFO, \objref{Table}{tab:5}); at 30\% adoption, the 350 kW scenario turns favorable (5.6\%; \objref{Table}{tab:5}). Normalizing by the 1.03 billion kWh of public charging that EPDK reports gives the same pattern (1.3\% and −7.6\% at 3\% adoption with reservation).

\begin{table}[!htb]
\centering
\caption{Case study results: four charger power levels at two plug-count configurations across Türkiye's 9,641-station EPDK network (full electrification, uncoordinated FIFO charging).}\label{tab:3}
\footnotesize
\setlength{\tabcolsep}{4pt}
\begin{tabular}{l rrrr rrrr}
\toprule
& \multicolumn{4}{c}{\textbf{50 plugs/station (482,050 total)}} & \multicolumn{4}{c}{\textbf{5 plugs/station (48,205 total)}} \\
\cmidrule(lr){2-5} \cmidrule(lr){6-9}
\textbf{Metric} & \textbf{5 kW} & \textbf{22 kW} & \textbf{50 kW} & \textbf{350 kW} & \textbf{5 kW} & \textbf{22 kW} & \textbf{50 kW} & \textbf{350 kW} \\
\midrule
\multicolumn{9}{l}{\emph{Charging operations}} \\
Avg.\ charge time/stop (min) & 497.2 & 113.1 & 49.8 & 7.2 & 497.2 & 113.1 & 49.8 & 7.2 \\
Avg.\ detour time/path (min) & 3.4 & 2.7 & 2.2 & 1.5 & 3.4 & 2.7 & 2.2 & 1.5 \\
Avg.\ stops/path & 2.4 & 2.4 & 2.4 & 2.5 & 2.4 & 2.4 & 2.4 & 2.5 \\
Additional distance (\%) & 0.59 & 0.47 & 0.39 & 0.26 & 0.59 & 0.47 & 0.39 & 0.26 \\
\multicolumn{9}{l}{\emph{Economic analysis (savings \%, positive = EV cheaper; waiting excluded)}} \\
Tier~1 (Operational) & 35.2 & 58.7 & 47.1 & 49.8 & 35.2 & 58.7 & 47.1 & 49.8 \\
Tier~2 (TCO) & 4.6 & 9.0 & 6.9 & 7.5 & 4.6 & 9.0 & 6.9 & 7.5 \\
Tier~3 (Infrastructure) & 4.6 & 9.0 & 6.8 & 7.0 & 4.6 & 9.0 & 6.8 & 7.0 \\
\multicolumn{9}{l}{\emph{Waiting and computation}} \\
FIFO mean wait/trip (days) & 261.8 & 61.7 & 28.2 & 4.5 & 2,635.3 & 630.5 & 294.5 & 55.2 \\
Runtime (min) & 1,323 & 1,308 & 1,330 & 1,310 & 1,325 & 1,303 & 1,319 & 1,341 \\
\bottomrule
\end{tabular}
\end{table}

\subsection{Ablation and Sensitivity Analysis}\label{sec:4_3}

To find what determines feasibility under a 20--80\% window, we ran an ablation on a 10\% sample of trips at 50 kW with about two DC plugs per station, starting from a restrictive configuration (trips must start within the window, SoC is tracked in 1\% buckets with each link's energy rounded up, and charging time is uncosted) and relaxing one choice at a time (\objref{Table}{tab:4}). The restrictive configuration makes only 0.3\% of trips feasible with a 20--80\% window but all of them with 10--100\%. Allowing a trip to start above the upper bound, as a driver who leaves home fully charged does, restores 100\% feasibility. Exact SoC tracking then cuts stops from 9.3 to 3.8 per trip, because rounding up overstates the energy of routes with many short links, and costing charging time lowers them to 3.6; more candidate stations change nothing. The low feasibility thus reflects the starting-SoC assumption, not gaps in the network.

\begin{table}[!htb]
\centering
\caption{Effect of three modeling choices on feasibility and charging behavior.}\label{tab:4}
\footnotesize
\begin{tabularx}{\textwidth}{>{\raggedright\arraybackslash}X cccccc}
\toprule
\textbf{Configuration} & \makecell{\textbf{SoC window}\\\textbf{(\%)}} & \makecell{\textbf{Feasible}\\\textbf{trips (\%)}} & \makecell{\textbf{Stops}\\\textbf{per trip}} & \makecell{\textbf{Charging time}\\\textbf{per stop (min)}} & \makecell{\textbf{Detour time}\\\textbf{per trip (min)}} \\
\midrule
(1) Restrictive configuration & 20--80 & 0.3 & 8.7 & 18.2 & 9.6 \\
(2) Restrictive configuration & 10--100 & 100.0 & 6.3 & 27.3 & 4.0 \\
(3) = (1) + start above the upper bound allowed & 20--80 & 100.0 & 9.3 & 17.4 & 7.7 \\
(4) = (3) + exact SoC tracking & 20--80 & 100.0 & 3.8 & 37.1 & 2.4 \\
(5) = (4) + charging time costed & 20--80 & 100.0 & 3.6 & 35.4 & 3.4 \\
(6) = (5) + up to three candidate stations (final model) & 20--80 & 100.0 & 3.6 & 35.4 & 3.4 \\
\bottomrule
\end{tabularx}

\smallskip
\raggedright 10\% sample of trips (224,753 that need charging under the 20--80\% window), 50 kW, about two DC plugs per station. Averages are over feasible trips.
\end{table}

Four univariate sweeps vary EV range (150--250, 250--400, 400--550, 550--700, and 700--900 km), the discharge rate upper bound (150, 160, 190, 220, and 250 Wh/km), charger power (50, 100, 250, 350, and 500 kW), and the operational SoC bounds (5--100\%, 10--100\%, 10--80\%, 20--100\%, and 20--80\%), holding the other parameters at their baseline values (\objref{Table}{tab:2}). Each level is run on a 1\% sample of trips at 350 kW with 50 plugs per station and 10 replications with distinct seeds, for 200 runs that typically take 20--25 minutes each.

The sensitivity analysis focuses on Tiers 1 and 2; Tier 3 depends on adoption rather than on the swept parameters and is reported in Sections~\ref{sec:4_2} and~\ref{sec:4_4}.

Discharge rate has the largest effect: Tier 1 savings rise from 49.8\% at 250 Wh/km to 58.9\% at 150 Wh/km (9.1 points) and Tier 2 savings from 7.5\% to 9.2\% (\objref{Figure}{fig:3} and \objref{Figure}{fig:4}), with no change in feasibility, stops, or detours (\objref{Figure}{fig:5}).

\begin{figure}[!htb]
\centering
\includegraphics[width=0.66\textwidth]{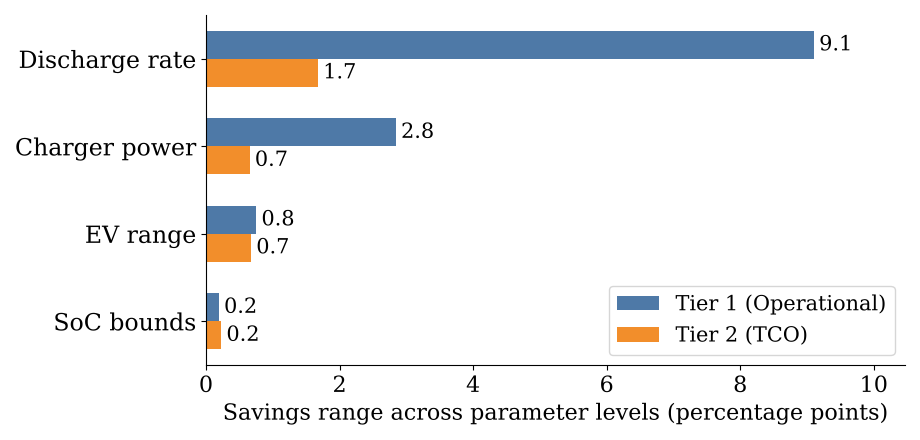}
\caption{Sensitivity range (in percentage points) of mean savings across parameter levels for Tier 1 (Operational) and Tier 2 (TCO).}\label{fig:3}
\end{figure}

EV range changes Tier 2 savings by only 0.7 points (7.0\% to 7.7\%); short-range vehicles make more stops (5.0 vs.~1.2 per trip) and longer detours (4.1 vs.~1.0 minutes), but all range levels are 100\% feasible.

\begin{figure}[!htb]
\centering
\begin{subfigure}[t]{0.48\textwidth}\centering\includegraphics[width=\linewidth]{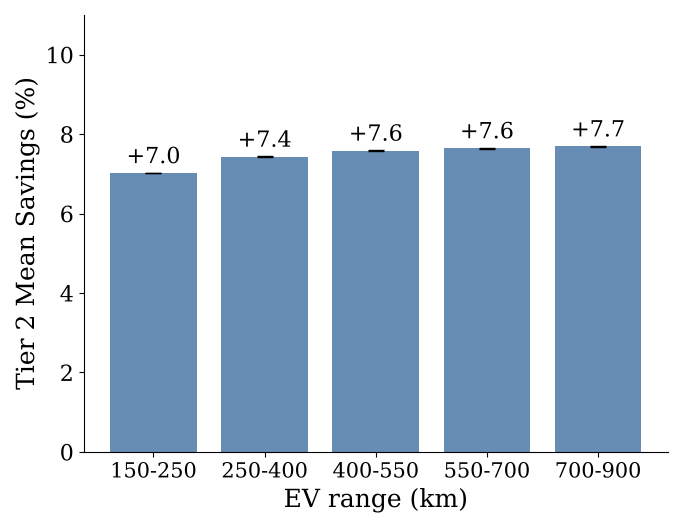}\caption{EV range}\end{subfigure}\hfill
\begin{subfigure}[t]{0.48\textwidth}\centering\includegraphics[width=\linewidth]{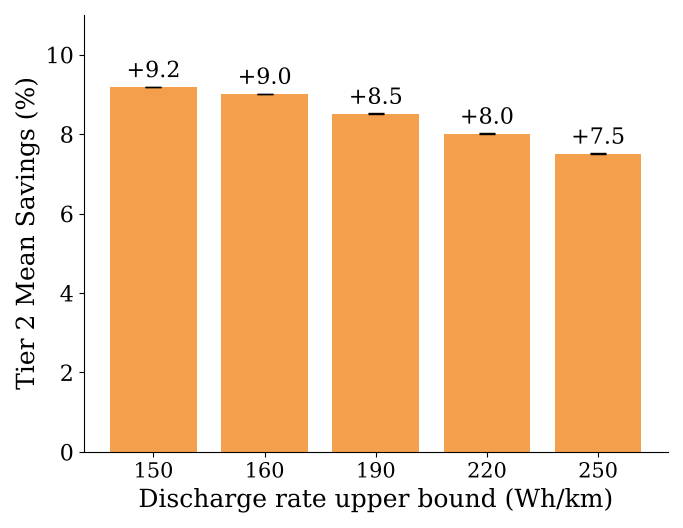}\caption{Discharge rate}\end{subfigure}

\medskip
\begin{subfigure}[t]{0.48\textwidth}\centering\includegraphics[width=\linewidth]{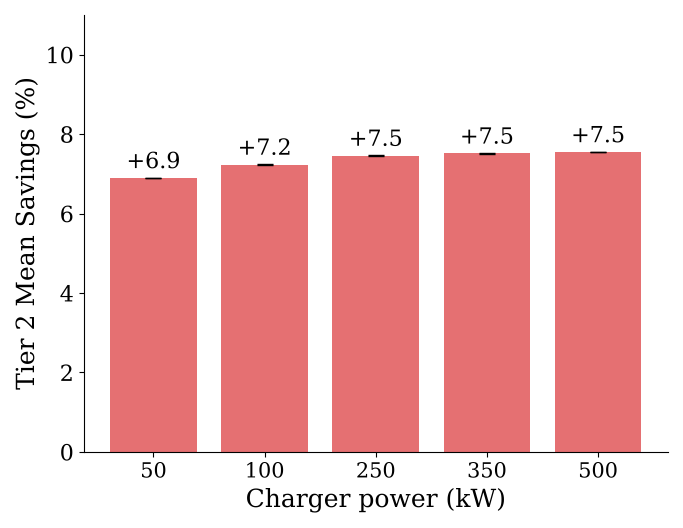}\caption{Charger power}\end{subfigure}\hfill
\begin{subfigure}[t]{0.48\textwidth}\centering\includegraphics[width=\linewidth]{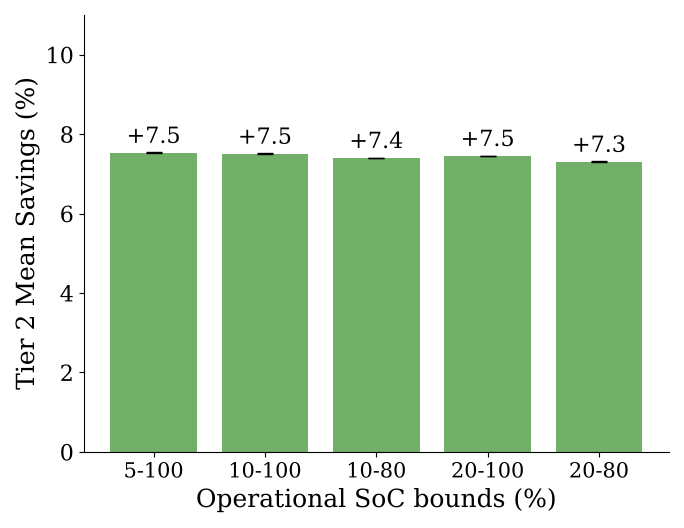}\caption{SoC bounds}\end{subfigure}
\caption{Tier 2 (TCO) mean savings by parameter level across the four sensitivity sweeps. Error bars show $\pm$1 standard deviation over 10 replications.}\label{fig:4}
\end{figure}

SoC bounds matter little once trips may start above the upper bound: all five windows, including 20--80\%, are 100\% feasible, and Tier 2 savings vary by 0.2 points. The 20--80\% window needs more stops (3.6 vs.~2.4 per trip at 5--100\%) and longer detours (2.8 vs.~1.6 minutes), which the network accommodates.

Charger power ranks second, changing Tier 1 savings by 2.8 points and Tier 2 savings by 0.7 points through the time cost of charging (49.7 minutes per stop at 50 kW against 5.1 at 500 kW). At a 1\% sample there is no congestion, so its larger role is as a capacity parameter (Section~\ref{sec:4_4}).

\begin{figure}[!htb]
\centering
\begin{subfigure}[t]{0.48\textwidth}\centering\includegraphics[width=\linewidth]{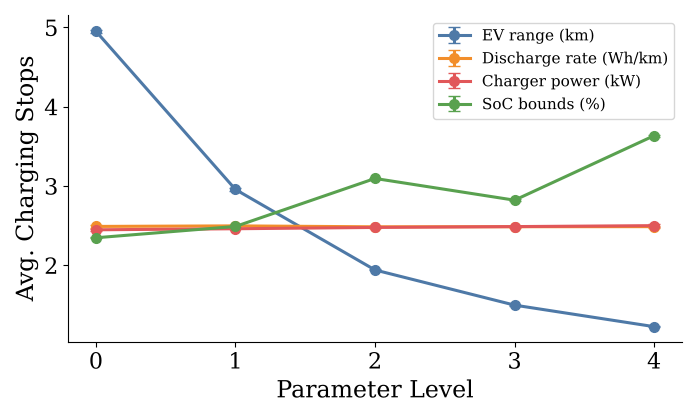}\caption{Average charging stops per path}\end{subfigure}\hfill
\begin{subfigure}[t]{0.48\textwidth}\centering\includegraphics[width=\linewidth]{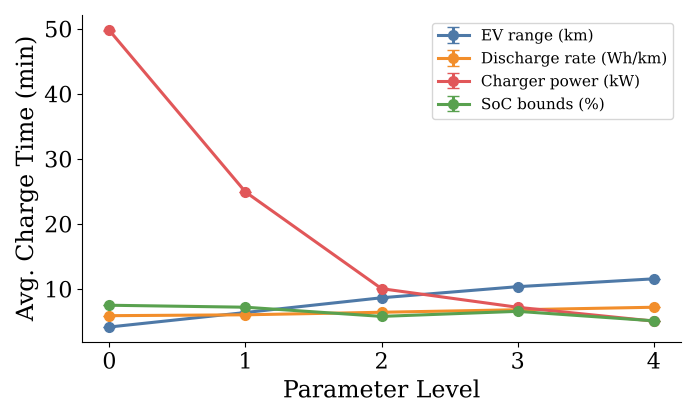}\caption{Average charge time per stop}\end{subfigure}

\medskip
\begin{subfigure}[t]{0.48\textwidth}\centering\includegraphics[width=\linewidth]{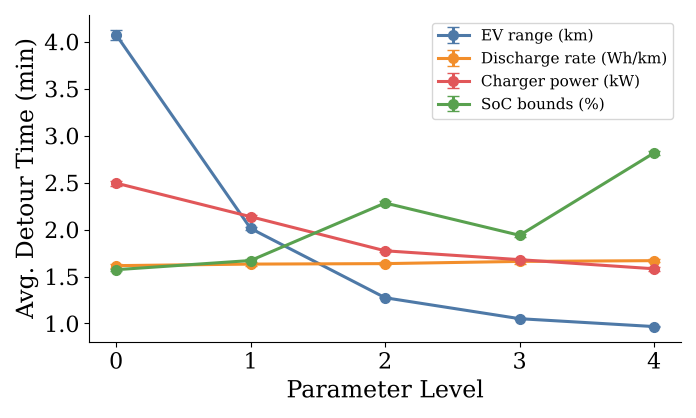}\caption{Average detour time per path}\end{subfigure}\hfill
\begin{subfigure}[t]{0.48\textwidth}\centering\includegraphics[width=\linewidth]{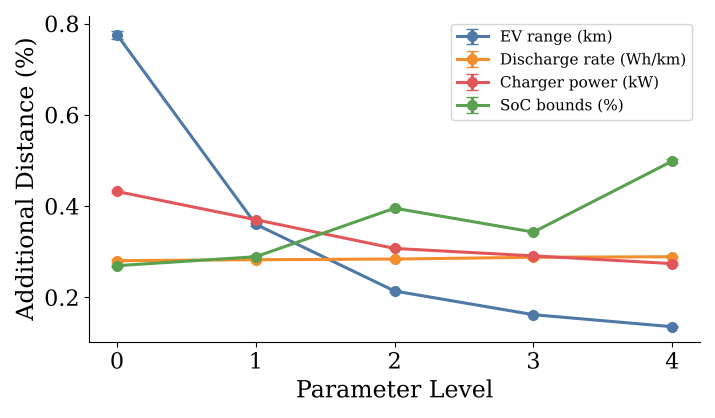}\caption{Additional travel distance}\end{subfigure}
\caption{Charging operations across sensitivity sweeps. Each line represents one parameter sweep; markers show mean values over 10 replications with $\pm$1 standard deviation error bars.}\label{fig:5}
\end{figure}

\subsection{Adoption and Congestion}\label{sec:4_4}

To examine congestion under realistic conditions, we set the plugs to the 19,663 DC sockets reported by EPDK (EPDK, \hyperlink{ref:epdk2026a}{2026a}), about two per station, and electrified 3\%, 10\%, 30\%, and 100\% of the trips that need charging; 3\% is close to today's EV share, and the first three levels correspond to about 6\%, 19\%, and 56\% of observed car vehicle-kilometers. Each level is run at 50 and 350 kW with FIFO charging (\objref{Table}{tab:5}). Sequential reservation is run where a run finishes within the 72-hour limit of a computing node: at 3\% adoption for both powers and at 10\% for 350 kW (\objref{Table}{tab:6}).

Uncoordinated charging produces long queues even at 3\% adoption: schedules converge on the same stations, mean waits reach 3.2 days at 350 kW and 20.8 days at 50 kW and grow with adoption, and 69--98\% of trips wait at least once. Reservation nearly eliminates waiting: at 3\% adoption, mean waits are 0.02 minutes at 350 kW and 1.3 minutes at 50 kW (95th percentiles of 0.04 and 2.7 minutes), and they stay below 0.1 minutes at 350 kW and 10\% adoption. The price is longer detours (0.37\% and 1.07\% of trip distance at 3\% adoption, against 0.28\% and 0.42\% under FIFO) and computing time. Because each vehicle is planned against all earlier reservations, the runtime grows much faster than the number of trips: at 350 kW, a reservation run takes about 8.5 hours at 3\% adoption and 36 hours at 10\%, against 48 minutes and 2.3 hours for FIFO, and at 30\% adoption (and at 10\% for 50 kW) it does not finish within 72 hours. Reservation results at higher adoption would therefore require a faster reservation algorithm or a decomposition of the national problem by region. Waiting does not enter the economic results, and where both runs exist, reservation changes the savings of FIFO by at most 0.7 percentage points (\objref{Table}{tab:5} and \objref{Table}{tab:6}).

\begin{table}[!htb]
\centering
\caption{Congestion and costs by EV adoption level under uncoordinated FIFO charging, with plug counts calibrated to EPDK (19,663 DC plugs, about two per station).}\label{tab:5}
\footnotesize
\begin{tabular}{rrrrrrrrr}
\toprule
\makecell{\textbf{Power}\\\textbf{(kW)}} & \makecell{\textbf{Adoption}\\\textbf{(\%)}} & \textbf{Trips} & \makecell{\textbf{Mean wait}\\\textbf{(days)}} & \makecell{\textbf{Trips}\\\textbf{waiting (\%)}} & \makecell{\textbf{Tier 1}\\\textbf{(\%)}} & \makecell{\textbf{Tier 2}\\\textbf{(\%)}} & \makecell{\textbf{Tier 3}\\\textbf{(\%)}} & \makecell{\textbf{Runtime}\\\textbf{(min)}} \\
\midrule
350 & 3 & 64,838 & 3.2 & 68.7 & 49.8 & 7.5 & $-$11.3 & 48 \\
350 & 10 & 216,124 & 12.9 & 81.4 & 49.8 & 7.5 & 1.9 & 138 \\
350 & 30 & 648,374 & 40.9 & 89.5 & 49.8 & 7.5 & 5.6 & 403 \\
350 & 100 & 2,161,247 & 140.3 & 94.8 & 49.8 & 7.5 & 7.0 & 1,356 \\
50 & 3 & 64,838 & 20.8 & 79.4 & 47.1 & 6.9 & 0.6 & 48 \\
50 & 10 & 216,124 & 71.8 & 89.4 & 47.1 & 6.9 & 5.0 & 138 \\
50 & 30 & 648,374 & 218.5 & 94.8 & 47.1 & 6.9 & 6.3 & 398 \\
50 & 100 & 2,161,247 & 738.4 & 98.0 & 47.1 & 6.9 & 6.8 & 1,328 \\
\bottomrule
\end{tabular}

\smallskip
\raggedright Trips: trips that need charging and are electrified. Mean wait per trip. Economic results (savings, positive = EV cheaper) exclude waiting.
\end{table}

\begin{table}[!htb]
\centering
\caption{Congestion and costs with sequential reservation, with plug counts calibrated to EPDK (19,663 DC plugs, about two per station).}\label{tab:6}
\footnotesize
\begin{tabular}{rrrrrrrrr}
\toprule
\makecell{\textbf{Power}\\\textbf{(kW)}} & \makecell{\textbf{Adoption}\\\textbf{(\%)}} & \textbf{Trips} & \makecell{\textbf{Mean wait}\\\textbf{(min)}} & \makecell{\textbf{95th percentile}\\\textbf{wait (min)}} & \makecell{\textbf{Tier 1}\\\textbf{(\%)}} & \makecell{\textbf{Tier 2}\\\textbf{(\%)}} & \makecell{\textbf{Tier 3}\\\textbf{(\%)}} & \makecell{\textbf{Runtime}\\\textbf{(min)}} \\
\midrule
350 & 3 & 64,838 & 0.02 & 0.04 & 49.8 & 7.4 & $-$11.4 & 520 \\
350 & 10 & 216,124 & 0.09 & 0.58 & 49.7 & 7.3 & 1.6 & 2,150 \\
50 & 3 & 64,838 & 1.3 & 2.7 & 46.6 & 6.3 & 0.0 & 510 \\
\bottomrule
\end{tabular}

\smallskip
\raggedright Same trips as the corresponding rows of \objref{Table}{tab:5}. Mean and 95th percentile wait per trip. Economic results (savings, positive = EV cheaper) exclude waiting.
\end{table}

\section{Conclusion}\label{sec:5}

This paper assessed whether Türkiye's 9,641 EPDK-licensed charging stations can support intercity EV travel by applying the scheduling framework of Cokyasar (\hyperlink{ref:cokyasar2026}{2026}) to 3 million synthetic trips on a 55,236-link highway network, through a case study of eight power and capacity scenarios, an adoption analysis with EPDK-calibrated plug counts, an ablation of three modeling choices, and a 200-run sensitivity analysis.

Four principal findings emerged. First, coverage is not the constraint: virtually all trips that need charging can be completed with the existing stations under both a 10--100\% and a 20--80\% SoC window, as long as a vehicle may start a trip above the upper bound of its window, and only 1.3\% of the network lies more than 50 km from a station. Second, capacity is the constraint: with about two DC plugs per station, uncoordinated charging builds long queues even at 3\% adoption, whereas reservation-based coordination keeps mean waits below two minutes. Third, EVs are about half as costly to operate as ICE vehicles at November 2026 prices (Tier 1 savings of 47--59\% at 22 kW and above), but their total cost of ownership advantage is narrow (4.6--9.0\%) and holds only while an EV costs at most about 6--11\% more than a comparable ICE vehicle. Fourth, the cost of charging infrastructure makes high-power charging uneconomic at today's adoption (−11.4\% at 350 kW and 3\% adoption) and favorable as adoption grows (5.6\% at 30\%).

For policy, the priority is capacity and coordination rather than coverage: adding plugs where trips converge and supporting booking services would help intercity travel more than new stations in most corridors, although the gaps in eastern and northeastern Anatolia (\objref{Figure}{fig:2}) still merit attention. The narrow ownership advantage makes the EV-to-ICE price gap, now favorable because of the lower special consumption tax on EVs, the main economic lever, and high-power chargers pay off only as adoption grows.

Several limitations should be noted. The scheduling framework assumes linear charging at a constant rate to the target SoC; nonlinear taper curves and partial charging would improve fidelity, and 350 kW should be read as an optimistic upper bound. Detours are approximated as Manhattan distances at a constant speed. Demand is distributed by population because no observed intercity origin--destination data exist; it overstates long-distance travel relative to observed traffic (about 1.9 times KGM's car vehicle-kilometers), so the congestion results are an upper bound, and income, car ownership, and tourism would shift trips between regions. Link speeds are free-flow speeds from OSM speed limits, with class defaults where limits are missing; actual speeds are lower, which affects travel and detour times but not energy use, which depends on distance. Plug counts are calibrated to national EPDK totals and spread evenly across stations because per-station counts are not published. The baseline excludes battery replacement, which common 8-year warranties cover for most of the vehicle's economic life, and all prices refer to a single date (November 2026). Finally, reservation represents an idealized booking service that all drivers use and honor, and its computing time limited it to adoption levels of up to 10\%.

Future work could replace the synthetic demand with observed origin--destination data, such as KGM traffic counts or mobile phone records, couple the scheduling framework with charger siting to target capacity additions, and extend the analysis to trucks and intercity buses.

\section*{Declarations}

\noindent\textbf{Conflict of interest.} The authors declare that there is no conflict of interest regarding the publication of this paper.

\noindent\textbf{Contributions.} Taner Cokyasar conducted study conception, design, data collection, analysis, and preparation. The author(s) reviewed and approved the manuscript.

\noindent\textbf{Research and publication ethics.} This study complies with research and publication ethics. It involves no human or animal participants and no personal data, so ethics committee approval or special permission was not required.

\noindent\textbf{Data and code availability.} The public inputs (OpenStreetMap, GADM, TÜİK, KGM, EPDK, and Şarj Gezgini) are cited with their sources and retrieval dates. The synthetic demand and the scheduling code rely on licensed POLARIS tools and cannot be redistributed under the license terms.

\noindent\textbf{Use of artificial intelligence.} Claude (Anthropic) was used to assist with code development, data processing, and language editing. The author reviewed and verified all text, code, results, and references and takes full responsibility for the content.

\section*{References}
\begin{refs}
\item\hypertarget{ref:aequilibraed2026}{}AequilibraE developers. (2026). \emph{AequilibraE: Python package for transportation modeling}. Retrieved November 16, 2026, from \url{https://www.aequilibrae.com}
\item\hypertarget{ref:auld2016}{}Auld, J., Hope, M., Ley, H., Sokolov, V., Xu, B., \& Zhang, K. (2016). POLARIS: Agent-based modeling framework development and implementation for integrated travel demand and network and operations simulations. \emph{Transportation Research Part C: Emerging Technologies}, \emph{64}, 101--116. \url{https://doi.org/10.1016/j.trc.2015.07.017}
\item\hypertarget{ref:bae2012}{}Bae, S., \& Kwasinski, A. (2012). Spatial and temporal model of electric vehicle charging demand. \emph{IEEE Transactions on Smart Grid}, \emph{3}(1), 394--403. \url{https://doi.org/10.1109/TSG.2011.2159278}
\item\hypertarget{ref:borlaug2020}{}Borlaug, B., Salisbury, S., Gerdes, M., \& Muratori, M. (2020). Levelized cost of charging electric vehicles in the United States. \emph{Joule}, \emph{4}(7), 1470--1485. \url{https://doi.org/10.1016/j.joule.2020.05.013}
\item\hypertarget{ref:chakraborty2019}{}Chakraborty, D., Bunch, D. S., Lee, J. H., \& Tal, G. (2019). Demand drivers for charging infrastructure---charging behavior of plug-in electric vehicle commuters. \emph{Transportation Research Part D: Transport and Environment}, \emph{76}, 255--272. \url{https://doi.org/10.1016/j.trd.2019.09.015}
\item\hypertarget{ref:cokyasar2026}{}Cokyasar, T. (2026). \emph{A national-scale EV charging scheduling framework: Optimal detour routing under infrastructure capacity constraints} (arXiv preprint arXiv:2607.18453). \url{https://arxiv.org/abs/2607.18453}
\item\hypertarget{ref:durmuenyapar2023}{}Durmuş Şenyapar, H. N., Akıl, M., \& Dokur, E. (2023). Adoption of electric vehicles: Purchase intentions and consumer behaviors research in Turkey. \emph{SAGE Open}, \emph{13}(2), 1--17. \url{https://doi.org/10.1177/21582440231180586}
\item\hypertarget{ref:epdk2026a}{}EPDK. (2026a). \emph{Şarj hizmeti piyasası aylık istatistikleri, Haziran 2026} {[}Charging service market monthly statistics, June 2026{]}. Enerji Piyasası Düzenleme Kurumu. Retrieved November 16, 2026, from \url{https://www.epdk.gov.tr}
\item\hypertarget{ref:epdk2026b}{}EPDK. (2026b). \emph{Elektrikli araç ve şarj altyapısı projeksiyonu} {[}Electric vehicle and charging infrastructure projection{]}. Enerji Piyasası Düzenleme Kurumu. Retrieved November 16, 2026, from \url{https://www.epdk.gov.tr}
\item\hypertarget{ref:evdatabase2025}{}EV Database. (2025). \emph{EV database: Energy consumption and range cheatsheet}. Retrieved from \url{https://ev-database.org/cheatsheet/energy-consumption-electric-car}
\item\hypertarget{ref:frade2011}{}Frade, I., Ribeiro, A., Gonçalves, G., \& Antunes, A. P. (2011). Optimal location of charging stations for electric vehicles in a neighborhood in Lisbon, Portugal. \emph{Transportation Research Record}, \emph{2252}(1), 91--98. \url{https://doi.org/10.3141/2252-12}
\item\hypertarget{ref:gadm2022}{}GADM. (2022). \emph{Database of Global Administrative Areas, version 4.1}. Retrieved July 10, 2026, from \url{https://gadm.org}
\item\hypertarget{ref:geofabrik2026}{}Geofabrik. (2026). \emph{OpenStreetMap data extract for Turkey}. Retrieved July 10, 2026, from \url{https://download.geofabrik.de/europe/turkey.html}
\item\hypertarget{ref:gnl2024}{}Gönül, Ö., Duman, A. C., \& Güler, Ö. (2024). A comprehensive framework for electric vehicle charging station siting along highways using weighted sum method. \emph{Renewable and Sustainable Energy Reviews}, \emph{199}, 114455. \url{https://doi.org/10.1016/j.rser.2024.114455}
\item\hypertarget{ref:gulbahar2023}{}Gulbahar, I. T., Sutcu, M., Almomany, A., \& Ibrahim, B. S. K. K. (2023). Optimizing electric vehicle charging station location on highways: A decision model for meeting intercity travel demand. \emph{Sustainability}, \emph{15}(24), 16716. \url{https://doi.org/10.3390/su152416716}
\item\hypertarget{ref:hagman2016}{}Hagman, J., Ritzén, S., Stier, J. J., \& Susilo, Y. (2016). Total cost of ownership and its potential implications for battery electric vehicle diffusion. \emph{Research in Transportation Business \& Management}, \emph{18}, 11--17. \url{https://doi.org/10.1016/j.rtbm.2016.01.003}
\item\hypertarget{ref:iea2024}{}IEA. (2024). \emph{Global EV outlook 2024}. Paris: International Energy Agency. Retrieved from International Energy Agency website: \url{https://www.iea.org/reports/global-ev-outlook-2024}
\item\hypertarget{ref:keil2017}{}Keil, P., \& Jossen, A. (2017). Calendar aging of NCA lithium-ion batteries investigated by differential voltage analysis and coulomb tracking. \emph{Journal of The Electrochemical Society}, \emph{164}(1), A6066--A6074. \url{https://doi.org/10.1149/2.0091701jes}
\item\hypertarget{ref:kgm2025}{}KGM. (2025). \emph{Karayolu ulaşım istatistikleri 2024} {[}Highway transportation statistics 2024{]}. Karayolları Genel Müdürlüğü. Retrieved November 16, 2026, from \url{https://www.kgm.gov.tr/SiteCollectionDocuments/KGMdocuments/Yayinlar/YayinPdf/KarayoluUlasimIstatistikleri2024.pdf}
\item\hypertarget{ref:knay2023}{}Kınay, Ö. B., Gzara, F., \& Alumur, S. A. (2023). Charging station location and sizing for electric vehicles under congestion. \emph{Transportation Science}, \emph{57}(6), 1433--1451. \url{https://doi.org/10.1287/trsc.2021.0494}
\item\hypertarget{ref:kuby2005}{}Kuby, M., \& Lim, S. (2005). The flow-refueling location problem for alternative-fuel vehicles. \emph{Socio-Economic Planning Sciences}, \emph{39}(2), 125--145. \url{https://doi.org/10.1016/j.seps.2004.03.001}
\item\hypertarget{ref:kullman2021}{}Kullman, N. D., Goodson, J. C., \& Mendoza, J. E. (2021). Electric vehicle routing with public charging stations. \emph{Transportation Science}, \emph{55}(3), 637--659. \url{https://doi.org/10.1287/trsc.2020.1018}
\item\hypertarget{ref:lee2019}{}Lee, Z. J., Li, T., \& Low, S. H. (2019). ACN-Data: Analysis and applications of an open EV charging dataset. \emph{Proceedings of the 10th ACM International Conference on Future Energy Systems (e-Energy '19)}, 139--149. \url{https://doi.org/10.1145/3307772.3328313}
\item\hypertarget{ref:li2016}{}Li, S., Huang, Y., \& Mason, S. J. (2016). A multi-period optimization model for the deployment of public electric vehicle charging stations on network. \emph{Transportation Research Part C: Emerging Technologies}, \emph{65}, 128--143. \url{https://doi.org/10.1016/j.trc.2016.01.014}
\item\hypertarget{ref:needell2016}{}Needell, Z. A., McNerney, J., Chang, M. T., \& Trancik, J. E. (2016). Potential for widespread electrification of personal vehicle travel in the United States. \emph{Nature Energy}, \emph{1}, 16112. \url{https://doi.org/10.1038/nenergy.2016.112}
\item\hypertarget{ref:nicholas2019}{}Nicholas, M. (2019). \emph{Estimating electric vehicle charging infrastructure costs across major U.S. metropolitan areas} (Working Paper 2019-14). International Council on Clean Transportation. \url{https://theicct.org/wp-content/uploads/2021/06/ICCT_EV_Charging_Cost_20190813.pdf}
\item\hypertarget{ref:numbeo2025}{}Numbeo. (2025). \emph{Cost of living in Türkiye}. Retrieved from \url{https://www.numbeo.com/cost-of-living/country_result.jsp?country=Turkey}
\item\hypertarget{ref:openstreetma2026}{}OpenStreetMap contributors. (2026). \emph{OpenStreetMap} {[}Data licensed under the Open Database License{]}. Retrieved July 10, 2026, from \url{https://www.openstreetmap.org}
\item\hypertarget{ref:petrolofisi2026}{}Petrol Ofisi. (2026). Akaryakıt fiyatları {[}Fuel prices{]}, İstanbul, November 1--15, 2026. Retrieved November 16, 2026, from \url{https://www.petrolofisi.com.tr/akaryakit-fiyatlari}
\item\hypertarget{ref:arjgezgini2026}{}Şarj Gezgini. (2026). \emph{Şarj Gezgini: Public charging station listing for Türkiye} {[}Third-party platform listing stations of EPDK-licensed charging network operators{]}. Data retrieved July 13, 2026, from \url{https://sarjgezgini.com/}
\item\hypertarget{ref:en2023}{}Şen, M., Yiğiter, M. S., \& Özcan, M. (2023). Why are consumers switching to electric vehicles? Analyzing consumers preferences for electric vehicles. \emph{Case Studies on Transport Policy}, \emph{14}, 101108. \url{https://doi.org/10.1016/j.cstp.2023.101108}
\item\hypertarget{ref:tang2017}{}Tang, W., \& Zhang, Y. J. (2017). A model predictive control approach for low-complexity electric vehicle charging scheduling: Optimality and scalability. \emph{IEEE Transactions on Power Systems}, \emph{32}(2), 1050--1063. \url{https://doi.org/10.1109/TPWRS.2016.2585202}
\item\hypertarget{ref:togg2026}{}Togg. (2026). T10X fiyat listesi, Kasım 2026 {[}T10X price list, November 2026{]}. Retrieved November 16, 2026, from \url{https://www.togg.com.tr}
\item\hypertarget{ref:toyotatrkiye2026}{}Toyota Türkiye. (2026). Corolla fiyat listesi, Kasım 2026 {[}Corolla price list, November 2026{]}. Retrieved November 16, 2026, from \url{https://www.toyota.com.tr}
\item\hypertarget{ref:trugo2026}{}Trugo. (2026). Şarj ücretleri {[}Charging prices{]}. Retrieved November 16, 2026, from \url{https://trugo.com.tr}
\item\hypertarget{ref:tk2024}{}TÜİK. (2024). \emph{Motorlu kara taşıtları İstatistikleri, aralık 2024}. Türkiye İstatistik Kurumu. Retrieved from Türkiye İstatistik Kurumu website: \url{https://data.tuik.gov.tr/Bulten/Index?p=Motorlu-Kara-Tasitlari-Istatistikleri-Aralik-2024-53762}
\item\hypertarget{ref:tk2025a}{}TÜİK. (2025a). \emph{Adrese dayalı nüfus kayıt sistemi sonuçları, 2024} {[}Address-based population registration system results, 2024{]} (Haber Bülteni No.~53783). Türkiye İstatistik Kurumu. Retrieved November 16, 2026, from \url{https://data.tuik.gov.tr/Bulten/Index?p=Adrese-Dayali-Nufus-Kayit-Sistemi-Sonuclari-2024-53783}
\item\hypertarget{ref:tk2025b}{}TÜİK. (2025b). \emph{Taşıt-kilometre istatistikleri, 2024} {[}Vehicle-kilometre statistics, 2024{]}. Türkiye İstatistik Kurumu. Retrieved November 16, 2026, from \url{https://data.tuik.gov.tr}
\item\hypertarget{ref:upchurch2009}{}Upchurch, C., Kuby, M., \& Lim, S. (2009). A model for location of capacitated alternative-fuel stations. \emph{Geographical Analysis}, \emph{41}(1), 85--106. \url{https://doi.org/10.1111/j.1538-4632.2009.00744.x}
\item\hypertarget{ref:witt2023}{}Witt, A. (2023). Determination of the number of required charging stations on a German motorway based on real traffic data and discrete event-based simulation. \emph{LOGI -- Scientific Journal on Transport and Logistics}, \emph{14}(1), 1--11. \url{https://doi.org/10.2478/logi-2023-0001}
\item\hypertarget{ref:zes2026}{}ZES. (2026). Fiyatlandırma {[}Pricing{]}. Retrieved November 16, 2026, from \url{https://zes.net/fiyatlandirma}
\item\hypertarget{ref:zhang2018}{}Zhang, T., Chen, X., Yu, Z., Zhu, X., \& Shi, D. (2018). A Monte Carlo simulation approach to evaluate service capacities of EV charging and battery swapping stations. \emph{IEEE Transactions on Industrial Informatics}, \emph{14}(9), 3914--3923. \url{https://doi.org/10.1109/TII.2018.2796498}
\end{refs}

\end{document}